\documentclass[11pt]{article}
\usepackage[cp1251]{inputenc}
\usepackage[ukrainian, english]{babel}

\usepackage{bezier,graphics,graphicx}
\usepackage{epstopdf}
\usepackage{amsmath,amsfonts,amssymb,amscd}
\usepackage{gensymb}

\begin{document}

\title{Directional diffusion splitting method for advection-diffusion-reaction model
\footnote{
advection-diffusion-reaction model, finite element method, Cauchy problem, parabolic equation, elliptic equation, boundary value problem, recurrent integration scheme, parallel algorithm, variational problem
}
}

\author{R. Drebotiy,~ H. Shynkarenko}

\date{Ivan Franko National University of Lviv,\\
1, Universytetska St., Lviv, 79000, Ukraine
\footnotetext{UDC 518:517.948}
}
\maketitle
\thispagestyle{empty}

\begin{abstract}
We propose certain approach of solving two-dimensional non-stationary and stationary advection-diffusion-reaction boundary value problems through their reduction to the set of corresponding one-dimensional problems. This method leverages special splitting and interpolation schemes, providing iterative algorithm with a large degree of parallelization possibilities. We combine that algorithm with the finite element method to solve obtained one-dimensional problems, but in fact, it can be combined also with other discretization methods, like finite volume or finite difference methods.
\end{abstract}

\section{Introduction}

In this article we consider two-dimensional non-stationary and stationary advection-diffusion-reaction (ADR) problems. There are several well-known general methods, which can be applied to the problems if such kind, to obtain numerical solution. For stationary case, we will have a boundary value problem (BVP) for equation of elliptic type, which can be successfully discretized by finite difference, finite volume and finite element methods. For the non-stationary problem we will have equation of parabolic type, combining boundary problem in space dimensions and the first-order initial problem in the time dimension. Typically such problems are discretized separately in time and space dimensions. According to the order of that discretization we can obtain Cauchy problem in standard form for the system of first-order equations or the sequence of elliptic problems in each time "slice". In space dimensions the finite elements or finite differences are used for discretization. In time dimension we typically use available toolset of first-order initial problem discretization methods, like Euler, Runge-Kutta, Adams methods etc.

For problems on large domains, modeling, for example, pollution migration over large country region, it is crucial to have efficient method of obtaining the solution, since the amount of discretization nodes and a size of resulting system of linear equation can be large. The same problem arises when we deal with singularly perturbed problems, having, for example, large advection-over-diffusion rate. In the last case many adaptive/stabilized schemes are proposed to deal with the singular perturbations. Without respect to which method/scheme is used, in general, for large problems, parallelization is the crucial tool, which is leveraged.

Not all methods are well suited for parallelization. Some which are suited, can have certain amount of synchronization points, so, in general, the level of parallelization capabilities can be different. For example well-known conjugate gradient (CG) method from scratch is not suited well to parallelization (as opposed to vector computations) \cite{CGpar} and so different techniques are implemented (like preconditioners) to increase the amount of parallelism and efficiency.

Also we may note, that iterative methods for large sparse systems, obtained from discretization, can often depend on a matrix condition number, which for singularly perturbed problems can be large, making iterative solvers less efficient.

Here we recall simple parallelizable method from \cite{DrebotiyCauchy} for pure advection-reaction problems and generalize it with the same ideas for the transport problem with diffusion present.

The paper is structured as follows: first we define model ADR problem; then we describe new algorithm; after that we present some numerical experiments.

\section{Advection-diffusion-reaction problem}
Let us consider the following two-dimensional Dirichlet problem for the stationary advection-diffusion-reaction equation:
\begin{equation}\label{BVP1}
\left\{ \begin{aligned}
   &\text{find function}~~u:{\bar\Omega}\to \mathbb{R}~~\text{such that:}  \\
   &-\mu \Delta u+\vec{\beta }\cdot \nabla u+\sigma u=f~~~\text{in}~~\Omega\subset\mathbb{R}^2,   \\
   &u=0~~~\text{on}~~\Gamma = \partial \Omega.  \\
\end{aligned} \right.
\end{equation}

Here $\Omega $ is a bounded domain with a Lipschitz boundary $\Gamma =\partial\Omega$, $\mu =const>0$ and $\sigma =const>0$ are coefficients of diffusion and reaction respectively, function $f=f(x)$ and vector $\vec\beta =(\beta_1(x), \beta_2(x))$ represent the sources and advection flow velocity respectively. We will consider non-compressible flow, i.e., $\nabla \cdot \vec{\beta }=0$ in $\Omega $. Also we consider the case, when vector field $\vec\beta$ does not have closed integral curves completely lying in $\bar\Omega$.

Problem \eqref{BVP1} is a limit case of corresponding non-stationary problem for parabolic equation:
\begin{equation}\label{Parabolic1}
\left\{ \begin{aligned}
   &\text{find function}~~u=u(x,t):{\bar\Omega \times [0,T] }\to \mathbb{R}~~\text{such that:}  \\
   &u_t^{\prime}-\mu \Delta_x u+\vec{\beta }\cdot \nabla_x u+\sigma u=f~~~\text{in}~~\Omega \times (0, T],   \\
   &u(x,t)=0,~~~(x, t)\in \Gamma \times [0, T],~~~\Gamma:=\partial\Omega \\
   &u(x,0)=u_0(x),~~~x\in \bar\Omega.  \\
\end{aligned} \right.
\end{equation}
Here we denoted corresponding operators by $x$ subscript to show that they are acting only in spatial dimensions. Also we let the time interval $[0, T]$ to be finite only for certainty. In general we can keep $T=+\infty$. For such case if we suppose, that the concentration $u$ stabilizes in time, i.e. we will have dynamic equilibrium, then $u_t^{\prime} \rightarrow 0$ as $t\rightarrow +\infty$, so we will have the degeneration of the problem \eqref{Parabolic1} to stationary problem \eqref{BVP1}.

In this article we propose the method for non-stationary problem \eqref{Parabolic1}. We can then use it iteratively through time integration scheme to obtain approximation for elliptic problem \eqref{BVP1}. For advection-dominated problem, we can specially choose the $u_0$ to make the time steps count needed for adequate approximation of the solution of stationary problem lower.

\section{Diffusion splitting}

Let us represent all vectors as columns by default. Recall, that the directional derivative of the function $u$ in the spatial direction $\vec e = (e_1, e_2)^T$, $\|\vec e\| = 1$ can be computed as $u_{\vec e}^{\prime} = \vec e \cdot \nabla_x u$. By the application of this formula two times, it is easy to see, that the second directional derivative in the same direction can be computed by the formula $u_{{\vec e} ^ 2}^{\prime\prime} = \nabla_x\cdot({\vec e}~{\vec e}^T\nabla_x u)$.

Let us consider the normalized advection vector field $\vec b (x) = (b_1(x), b_2(x)) := \vec\beta / \|\vec\beta\|$. Consider also orthogonal vector field $\vec \gamma(x) = (\beta_2(x), -\beta_1(x))$ and corresponding normalized field $\vec p = \vec \gamma / \|\vec \gamma\|$.

It is obvious, that identity matrix $I$ can be then decomposed as a sum of two orthogonal projections:
\begin{equation}\label{IdentDecomp}
I = {\vec b}~{\vec b}^T + {\vec p}~{\vec p}^T.
\end{equation}
Now we can transform the first (diffusion) term of the ADR equation:
\begin{equation}\label{DiffTerm}
\mu \Delta_x u = \mu \nabla_x \cdot (I\nabla_x u) = \mu \nabla_x \cdot ({\vec b}~{\vec b}^T\nabla_x u) + \mu \nabla_x \cdot ({\vec p}~{\vec p}^T\nabla_x u)
\end{equation}
Note, that each of obtained terms represents second directional derivative of $u$ in the appropriate direction.

\section{Variational formulation and semi-discretization in time}

As a backend of full discretization of our problems we use finite element method and thus we need to reformulate our problem in the form of variational equation.

The boundary value problem \eqref{BVP1} admits the following variational formulation:
\begin{equation}\label{VP1}
\left\{ \begin{aligned}
   &\text{find}~~u\in V:=H_{0}^{1}(\Omega )~~\text{such that},  \\
   &a(u,v)=\langle l,v\rangle ~~~\forall v\in V,  \\
\end{aligned} \right.
\end{equation}
where:
\begin{equation}\label{VP1a}
\left\{ \begin{aligned}
   &a(u,v)=\int\limits_{\Omega }{(\mu \nabla u\cdot \nabla v+\vec{\beta }v\cdot \nabla u+\sigma uv)}dx~~~\forall u,v\in V,\\
   &\langle l,v\rangle =\int\limits_{\Omega }{f}vdx~~~\forall v\in V.\\
\end{aligned} \right.
\end{equation}
Let us define time step $\Delta t$ and a parameter $\theta \in (0, 1)$. Consider standard Lebesque scalar product $(w, q) := \int_\Omega wqdx$ Using components \eqref{VP1a}, and approach from \cite{Trush} we can formulate now one-step recurrent scheme for semi-discretization in time of non-stationary problem \eqref{Parabolic1}:
\begin{equation}\label{SemiDiscr}
\left\{ \begin{aligned}
   &(\dot u_{j+\frac{1}{2}}, v) + \theta \Delta t a(\dot u_{j+\frac{1}{2}}, v) = \langle l,v\rangle - a(u_j, v) \\
   &u_{j+1} = u_j + \Delta t \dot u_{j+\frac{1}{2}}, ~~~j=0,1,...,\\
\end{aligned} \right.
\end{equation}
where $u_j$ is an approximation to the function $u(x, t_j)$, $t_j = j\Delta t$.

\section{Semi-discretization splitting}

Let us substitute formula \eqref{DiffTerm} into bilinear form $a$. Let us rewrite it in the following way:
\begin{equation}\label{VPcompsplit1}
a(u,v):=s(u,v)+m(u,v),
\end{equation}
where
\begin{equation}\label{VPcompsplit2}
\left\{ \begin{aligned}
   &s(u,v):=\int\limits_{\Omega }{(\mu ({\vec b}^T\nabla_x u)({\vec b}^T\nabla_x v)+v\|\vec{\beta }\|{\vec b}^T\nabla_x u+\sigma uv)}dx\\
   &m(u,v):=\int\limits_{\Omega }{\mu ({\vec p}^T\nabla_x u)({\vec p}^T\nabla_x v)}dx\\
\end{aligned} \right.
\end{equation}
We propose to use instead of recurrent scheme \eqref{SemiDiscr} the following two-step scheme:
\begin{equation}\label{SemiDiscrSplit}
\left\{ \begin{aligned}
   &(\dot u_{j+\frac{1}{4}}, v) + \theta \Delta t s(\dot u_{j+\frac{1}{4}}, v) = \langle l,v\rangle - s(u_j, v) \\
   &u_{j+\frac{1}{2}} = u_j + \Delta t \dot u_{j+\frac{1}{4}},\\
   &(\dot u_{j+\frac{3}{4}}, v) + \theta \Delta t m(\dot u_{j+\frac{3}{4}}, v) = - m(u_{j+\frac{1}{2}}, v) \\
   &u_{j+1} = u_{j+\frac{1}{2}} + \Delta t \dot u_{j+\frac{3}{4}}, ~~~\forall v\in V,~~~j=0,1,...,\\
\end{aligned} \right.
\end{equation}
The general idea behind that splitting is to decouple diffusion component, that is orthogonal to the direction of advection. We will show in the next section, that such splitting lead us to effective way of solving obtained variational problems in parallel. The idea of \textit{operator splitting} is known in the literature. For example we can recall \textit{Chorin's projection method} for Navier-Stokes equations where computations of the velocity and the pressure fields are decoupled. Proposed approach of splitting in its root is similar to general splitting approach based on the application of Baker–Campbell–Hausdorff formula, which gives us the possibility of switching from one time step into two steps with individual terms of original operator (with asymptotically quadratic error in time step length).

Note, that if we wish to solve original stationary problem with dominated advection, we should consider corresponding non-stationary counterpart with $u_0 \in V$ calculated as a solution of the following equation:
\begin{equation}\label{u0stationary}
s(u_0, v) = \langle l,v\rangle ~~~\forall v\in V
\end{equation}
since in that way we will obtain "from start" good approximation to the initial equation, since the dynamic of the processes will be mainly directed by advection in that case.

\section{One-dimensional reduction of semi-discretization substeps}

We denote by $\vec n (x)$ the vector of the unit outward normal to the boundary $\partial \Omega$ in point $x$. Let us consider the sets:
\begin{equation}\label{cauchyStartSets}
\begin{aligned}
   &\Gamma_{\vec\beta}:=\{x\in\partial\Omega|\vec n (x) \cdot \vec\beta(x) < 0\}\\
   &\Gamma_{\vec\gamma}:=\{x\in\partial\Omega|\vec n (x) \cdot \vec\gamma(x) < 0\}\\
\end{aligned}
\end{equation}
Consider the functions $x=x(t, \xi)$ and $y=y(t, \eta)$, $x, y \in \Omega, t\ge0, \xi,\eta\in[0, 1]$, such that $x(0, \xi)=\rho(\xi)\in\Gamma_{\vec\beta}$ and $y(0, \eta)=\tau(\xi)\in\Gamma_{\vec\gamma}$ are parametrizations of curves $\Gamma_{\vec\beta}$ and $\Gamma_{\vec\gamma}$ correspondingly.
Consider the Cauchy problems for those functions:
\begin{equation}\label{CauchyBeta}
\left\{ \begin{aligned}
   &x_t^{\prime} = \vec\beta(x),\\
   &x(0, \xi)=\rho(\xi)\\
\end{aligned} \right.
\end{equation}
and
\begin{equation}\label{CauchyGamma}
\left\{ \begin{aligned}
   &y_t^{\prime} = \vec\gamma(y),\\
   &y(0, \eta)=\tau(\xi).\\
\end{aligned} \right.
\end{equation}
Those problems define the families of integral curves of $\vec\beta$ and $\vec\gamma$ covering the domain $\Omega$.
Consider now some finite sets of those curves $B=\{x(t, \xi_i)\}_{i=1}^n$ and $G=\{y(t, \eta_i)\}_{i=1}^m$

We now reduce the substeps of \eqref{SemiDiscrSplit} to those sets of curves.
Let use rewrite the forms from \eqref{VPcompsplit2} using the mentioned earlier notion of directional derivative:
\begin{equation}\label{VPcompsplit3}
\left\{ \begin{aligned}
   &s(u,v):=\int\limits_{\Omega }{(\mu u_{\vec b}^{\prime} v_{\vec b}^{\prime}+ \|\vec{\beta }\|u_{\vec b}^{\prime}v+\sigma uv)}dx\\
   &m(u,v):=\int\limits_{\Omega }{\mu u_{\vec p}^{\prime} v_{\vec p}^{\prime}}dx\\
\end{aligned} \right.
\end{equation}

Consider for example some curve $L_i\in B$. Let us set $v(x) = w_i(x) \delta(dist(x, L_i))$, where $w_i\in H_0^1([0, |L_i|])$ and $\delta$ is a Dirac $\delta$-function. Taking into account the integration by parts formula and the fact, that Dirac $\delta$-function can be expressed as weakly convergent sequence of smooth functions, we can express the first two equations from \eqref{SemiDiscrSplit} for the curve $L_i$ in the following way:
\begin{equation}\label{SemiDiscrSplit1Dbeta}
\left\{ \begin{aligned}
   &(\dot u_{j+\frac{1}{4}}^i, w_i) + \theta \Delta t s_i(\dot u_{j+\frac{1}{4}}^i, w_i) = \langle l_i, w_i\rangle - s_i(u_j^i, w_i) \\
   &u_{j+\frac{1}{2}}^i = u_j^i + \Delta t \dot u_{j+\frac{1}{4}}^i, ~~~\forall w_i\in H_0^1([0, |L_i|]),~~~j=0,1,..., i=1..n\\
\end{aligned} \right.
\end{equation}
where
\begin{equation}\label{SemiDiscrSplit1Dbeta1}
\left\{ \begin{aligned}
   &s_i(u,v):=\int\limits_{L_i}{(\mu u^{\prime} v^{\prime}+ \|\vec{\beta }\|u^{\prime}v+\sigma uv)}dl\\
   &\langle l_i,v\rangle =\int\limits_{L_i }{f}vdl~~~\forall u,v\in H_0^1([0, |L_i|]).\\
\end{aligned} \right.
\end{equation}

In the same way we can consider arbitrary curve $K_i$. Let us set $v(x) = w_i(x) \delta(dist(x, K_i))$, where $w_i\in H_0^1([0, |K_i|])$. Now, we can express the second two equations from \eqref{SemiDiscrSplit} for the curve $K_i$ in the following way:

\begin{equation}\label{SemiDiscrSplit1Dgamma}
\left\{ \begin{aligned}
   &(\dot u_{j+\frac{3}{4}}^i, w_i) + \theta \Delta t m_i(\dot u_{j+\frac{3}{4}}^i, w_i) = - m_i(u_{j+\frac{1}{2}}^i, w_i) \\
   &u_{j+1}^i = u_{j+\frac{1}{2}}^i + \Delta t \dot u_{j+\frac{3}{4}}^i, ~~~\forall w_i\in H_0^1([0, |K_i|]),~~~j=0,1,..., i=1..m\\
\end{aligned} \right.
\end{equation}
where
\begin{equation}\label{SemiDiscrSplit1Dgamma1}
m_i(u,v):=\int\limits_{K_i}{\mu u^{\prime} v^{\prime}}dl.
\end{equation}

In the \eqref{SemiDiscrSplit1Dbeta} and \eqref{SemiDiscrSplit1Dgamma} all functions with the upper index $i$ are exactly equal on corresponding curves to the appropriate functions without that index.

To solve approximately the obtained 1D problem we can use finite element method.

So, we transformed our two steps from \eqref{SemiDiscrSplit} into two sets of BVPs on the sets of orthogonal curves. In that way we need to have a way to interpolate the values from one set to another one, to be able to substitute the values from $u_{j+\frac{1}{2}}$ into the second step in \eqref{SemiDiscrSplit}.
In ideal case it will be good to have the corresponding intersection points of those orthogonal curves, but in practice finding such points is not so effective. We propose here the similar way as in \cite{DrebotiyCauchy}, which we recall here.

Let us define bounding box for domain   as:
\begin{equation}\label{BndBox}
	P=[\underset{(x,y)\in \Omega }{\mathop{\min }}\,x,\,\,\underset{(x,y)\in \Omega }{\mathop{\max }}\,x]\times [\underset{(x,y)\in \Omega }{\mathop{\min }}\,y,\,\,\underset{(x,y)\in \Omega }{\mathop{\max }}\,y]
\end{equation}
Let us cover this bounding box by the rectangular mesh of cells, defined by $M$ horizontal and $N$ vertical lines.
To obtain the value in grid point we do the following. First, using elementary geometry, we find points of intersections between out integral curve segment and the grid lines. We can keep in memory data structure with lines and their segments with all needed data and update that structure on the fly, when we found the next FE approximation on certain integral curve. For each point we found, we use linear interpolation to compute interpolated value of corresponding function (for example $u$) in that point. For each line segment we can keep two pair of values. First pair will correspond to intersection point with minimal coordinate with appropriate function value. The second pair will correspond to maximum intersection coordinate on that segment.
After we found all FE approximations to 1D problems, we will have populated arrays of lines and their segments. Each segment will contain min/max intersection points. Now we can identify for the point of the grid two closest points on corresponding horizontal line and two on a vertical line. They will be corresponding min/max points from the segment data structures.
Now we can use linear interpolation on horizontal line to interpolate value of function in the grid point using closest point on that line. The same thing we can do on vertical line. Now to obtain final result we just take average two values which we found and treat that average as an appropriate approximation to the solution in the grid node.

\subsection{Bilinear interpolation}
Now, when we interpolated values into rectangular grid, we need to have possibility to recover the value in arbitrary point, to be able to map values generated on $\vec\beta$ integral curves onto appropriate $\vec\gamma$ curves and vice versa. For that we can use bilinear interpolation on each grid cell to find approximate values in all other points.
Suppose we have point $(x,y)\in \Omega $, corresponding bounding box $P=[a,b]\times [c,d]$ and we have grid with ${{k}_{x}}$ and ${{k}_{y}}$ cells in row and column respectively. If we enumerate cells along the coordinate axes using two indexes $({{i}_{x}},{{i}_{y}})$, then we can find trivially the cell indexes for the point:
\begin{equation}\label{bilinInd}
{{i}_{x}}=\left\lceil \frac{x-a}{b-a}{{k}_{x}} \right\rceil ,\,\,\,\,\,\,\,\,\,\,\,\,\,\,\,\,\,\,{{i}_{y}}=\left\lceil \frac{x-c}{d-c}{{k}_{y}} \right\rceil
\end{equation}
Suppose, that we have the following corner values of ${{u}_{0}}$ for cell: ${{p}_{11}},\,\,{{p}_{12}},\,\,{{p}_{21}},\,{{p}_{22}}$ (bottom left, bottom right, top left, top right). Then, we can interpolate those values to find approximate value of target function in the point $(x,y)$:
\begin{equation}\label{bilinInterp}
I(x,y)=(1-{{\lambda }_{y}})(1-{{\lambda }_{x}}){{p}_{11}}+(1-{{\lambda }_{y}}){{\lambda }_{x}}{{p}_{12}}+{{\lambda }_{y}}(1-{{\lambda }_{x}}){{p}_{21}}+{{\lambda }_{y}}{{\lambda }_{x}}{{p}_{22}},
\end{equation}
where:
\begin{equation}\label{bilinInterpComp}
{{\lambda }_{x}}=\left\{ \frac{x-a}{b-a}{{k}_{x}} \right\},\,\,\,\,\,\,\,\,\,\,\,\,\,{{\lambda }_{y}}=\left\{ \frac{x-c}{d-c}{{k}_{y}} \right\}
\end{equation}
and by $\{\}$ we denote fractional part of a number.

\section{Numerical experiment}

Let us consider 2D stationary problem with the following data:
\begin{equation}\label{Data2D}
\Omega=(0,1)^2, \mu=1, \vec{\beta}(x,y)=(-5(y+1), 5(x+1))^T, \sigma=1, f(x,y)=5
\end{equation}
We used parameter $\theta=0.5$, which, in fact, defines classic Crank–Nicolson scheme.

\begin{figure}[h!]
    \centering
    \includegraphics[width=0.55\linewidth]{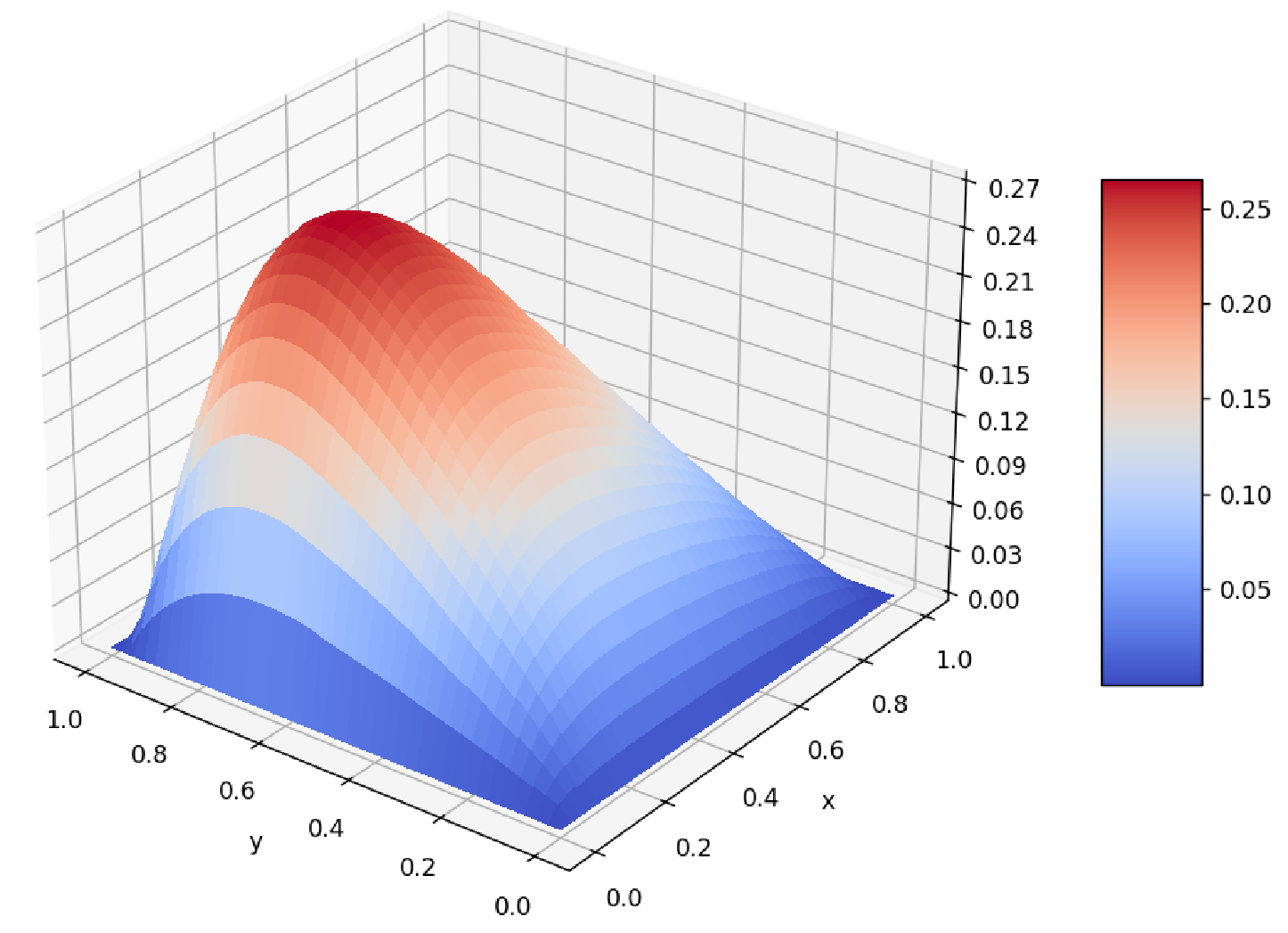}%
    \caption{Approximate solution of problem with data \ref{Data2D}, obtained by the proposed method.}
    \label{numFig}
\end{figure}

We solved this problem by using corresponding non-stationary problem with $u_0$ selected as a solution of equation \eqref{u0stationary} and executing 50 time iterations with $\Delta t = 0.001$.
Obtained solution was coinciding with the one, obtained using finite element method with FEniCS library on the uniform triangular mesh with 450 elements. We calculated relative error with respect to the solution by FEniCS library in $L^\infty$ and $L^1$ norms, which was 5\% and 2\% respectively. Plot of obtained approximation is depicted on Fig. \ref{numFig}.

Note, that proposed method can be implemented with the great level of parallelism with small number synchronization points. We can obviously solve in parallel all obtained 1D problems and, as mentioned in \cite{DrebotiyCauchy}, such interpolation scheme with intermediate grid can be implemented also in parallel way.

\section{Conclusions}

In this article we constructed splitting method for advection-diffusion-reaction model, which can be used to transform the two-dimensional problem to the set of one-dimensional problems and solve them with a high level of parallelism.

Note, that it looks like, that for singularly perturbed problems, we can simply use in parallel any adaptive or stabilization scheme for each one of obtained one-dimensional problems along advection field. This approach was not implemented yet and it is considered for further research.

Also one of the improvements of the constructed scheme is to get rid of that bilinear interpolation part by mapping somehow intersection nodes between two sets of integral curves. It seems reasonable to leverage the fact, that in each intersection point the segments of lines are orthogonal. One idea is to construct the set of "buckets" containing groups of segments with close "azimuthal angles". In that way we can quickly find the segments of the opposite curve group which are intersecting with the given segment (from the other set of orthogonal curves). This idea is also not implemented yet and it is considered for further research/improvement.


\begin{thebibliography}{99}

\bibitem{CGpar}
    \emph{Dianne P. O'LEARY}
    Parallel implementation of the block conjugate gradient algorithm // Parallel Computing 5 (1987) 127-139
\bibitem{DrebotiyCauchy}
    \emph{Drebotiy\,R., Shynkarenko\,H.}
    Fully parallel algorithm for solution of advection-reaction Cauchy problem in one finite element regularization scheme // Manufacturing Processes. Actual Problems.– Politechnika Opolska.– Opole, 2023.– Vol. 1.– Pp. 35-44.
\bibitem{Trush}
    \emph{Trushevsky\,V.\,M., Shynkarenko\,H.\,A., Shcherbyna\,N.\,M.}
    Finite element method and artificial neural networks: theoretical aspects and application. //
    Lviv: Ivan Franko National University of Lviv, 2014, ISBN 978-617-10-0127-5 (in Ukrainian)
\end{thebibliography}
\end{document}